\documentclass{amsart}

\usepackage{xcolor}
\usepackage{amssymb}
\usepackage{latexsym}
\usepackage{amsmath}
\usepackage{euscript}
\usepackage{graphics}
\usepackage{epigraph}

   \def\dZ{{\mathbb Z}}

\def\bm\chi{\mbox{\boldmath$\chi$}}

\let\xker=\ker \def\ker{{\xker\,}}

\def\Im{\operatorname{Im}}

\newtheorem{theorem}{Theorem}[section]
\newtheorem{proposition}[theorem]{Proposition}

\theoremstyle{remark}

\newtheorem{remark}[theorem]{Remark}

\numberwithin{equation}{section}

\author[M.~Derevyagin]{Maxim~Derevyagin}
\address{
MD,
Department of Mathematics\\
University of Connecticut\\
341 Mansfield Road, U-1009\\
Storrs, CT 06269-1009, USA}
\email{maksym.derevyagin@uconn.edu}

\date{\today}

\dedicatory{\it Dedicated with gratitude to my teacher and advisor, Vladimir A. Derkach, \\on the occasion of his 75th birthday. }

\begin{document}

\title[Thiele Meets Schur]{Reflections on Thiele, Schur, and Continued Fractions}

\begin{abstract} 
We survey several classical continued-fraction constructions from a common point of view centered on linear-fractional transformations, interpolation, and the recurrence relations they generate. Beginning with Thiele's interpolation algorithm, we trace how similar mechanisms reappear in the Schur algorithm, Wall's continued fractions, and multipoint Nevanlinna--Pick interpolation. Particular attention is paid to the emergence of \(R_{II}\)-type continued fractions, which arise naturally as even parts of Thiele continued fractions, and to the associated pencils of tridiagonal matrices. A detour through Euler's continued-fraction formula illustrates how rational Nevanlinna functions lead naturally to finite \(R_{II}\)-pencils, while the Schur algorithm is related to Wall pencils and, through a different reorganization, to CMV matrices. We also discuss a modified Schur transformation that removes an interpolation condition at a pair of symmetric points, which leads to a linear pencil of two Hermitian tridiagonal matrices. Finally, we indicate how periodically occurring interpolation nodes at infinity convert \(R_{II}\)-pencils into ordinary banded matrices, suggesting a connection with GMP matrices.
\end{abstract}
\maketitle

\section*{Prologue}
\addcontentsline{toc}{section}{Prologue}
{\it At this point I no longer remember whether I first encountered the Nevanlinna--Pick problem or the moment problem. I certainly learned about continued fractions before either of them, and it always seemed somewhat unsatisfactory to me that there was no analogue of the Jacobi continued fraction for Nevanlinna--Pick interpolation, comparable to the role it plays in the theory of moment problems. Eventually, in the course of working on these questions, I came to understand what such an analogue should look like, learned more of the history surrounding it, and realized that Wall had already come surprisingly close to this point of view.

This survey presents a largely historical and algebraic perspective on that journey. Behind the classical interpolation problems discussed here lies a spectral picture more general than that of the moment problem: one is naturally led to operator pencils, Hermitian operators, and their corresponding Weyl functions. The single Jacobi operator and its Weyl $m$-function are recovered when the Nevanlinna–Pick problem reduces to a moment problem. The history of CMV matrices provides a useful reminder that even a very classical theory may conceal a particularly natural matrix representation for a long time. As one of the detours in this survey shows, CMV matrices themselves can be brought into the same circle of continued-fraction and Schur-type constructions. In this sense, it is perhaps not surprising that some of the analogous structures associated with Nevanlinna--Pick interpolation remained less visible, even though Wall had already obtained a closely related construction in a setting not very far removed from the classical moment problem.}

\section{Introduction}

Continued fractions have played several distinct roles in the development of analysis. Emerging from the Euclidean algorithm and the study of rational approximations to irrational numbers, they were developed systematically by Euler and subsequently used by Lambert, Lagrange, and Gauss in problems ranging from Diophantine approximation to the representation of analytic and special functions. During the 19th century, the work of Jacobi, Chebyshev, Markov, and Stieltjes revealed increasingly deep connections between continued fractions, recurrence relations, moment problems, and orthogonal polynomials. From this point of view, continued fractions may be regarded as a natural mechanism for encoding sequences of rational approximants and the recurrences they satisfy. For example, Newton's divided-difference formula provided a systematic way to pass from interpolation data to a polynomial or Newton-series expansion, while Euler showed that suitable series expansions could be transformed into continued fractions. Thiele subsequently introduced a direct interpolation calculus for continued fractions: replacing divided differences by reciprocal differences, he constructed continued fractions directly from interpolation data and thereby obtained a rational counterpart of Newton interpolation.
Shortly after that, the work of Schur and Nevanlinna brought an additional function-theoretic principle into this picture: interpolation conditions are removed recursively by linear-fractional transformations chosen so that a distinguished class of analytic functions is preserved. Various historical accounts of continued fractions and their connections with approximation theory, orthogonal polynomials, and analytic function theory can be found, for example, in the books of Brezinski \cite{Brezinski} and Khrushchev \cite{Khrushchev08}.

The Schur and Nevanlinna algorithms, however, are not continued-fraction algorithms in their original form: their defining feature is a recursive linear-fractional transformation that preserves a distinguished class of analytic functions. Nevertheless, as was already implicit in Wall's treatment and becomes explicit after a suitable reformulation, the Schur and Nevanlinna algorithms can be viewed as linear-fractional transformations that generate continued fractions. Thus the continued-fraction and Schur--Nevanlinna viewpoints are not merely analogous but, in an appropriate form, equivalent. Thiele's interpolation scheme provides a natural setting in which to make this connection explicit. This motivates returning to Thiele's interpolation fraction and asking how its interpolation mechanism interacts with such class-preserving transformations.

We begin with the interpolation problem underlying Thiele's construction. Let

$$
z_0,z_1,\ldots,z_N, \ldots
$$
be distinct interpolation nodes and let
$$
w_0,w_1,\ldots,w_N, \ldots
$$
be prescribed values.  We seek a rational function \(f\) satisfying
\begin{equation}\label{interp}
f(z_j)=w_j,\qquad j\in\dZ_+.
\end{equation}
Under the usual nondegeneracy assumptions, Thiele's algorithm constructs such a function in the form
\begin{equation}\label{thiele}
f(z)=
a_0+
\cfrac{z-z_0}{
a_1+
\cfrac{z-z_1}{
a_2+
\cfrac{z-z_2}{
\ddots+
\cfrac{z-z_{N-1}}{a_N+\cfrac{z-z_N}{
\ddots}}}}},
\end{equation}
where the numbers \(a_j\) are the reciprocal differences associated with the interpolation data. Evidently, \eqref{thiele} is, in general, only a formal continued fraction, and the existence of an interpolated function naturally raises the question of its convergence. For a finite interpolation set, however, the construction produces a rational solution of the interpolation problem and, upon allowing the terminal term to vary, generates a family of solutions.

Conversely, every solution to a problem with a finite interpolation set can be reduced by reversing the construction, showing that this family in fact exhausts the solution set. Namely, let us put

$$
f_0=f,\qquad a_0=f_0(z_0),
$$
and define
\begin{equation}\label{thiele-reduction}
f_{k+1}(z)
=\frac{z-z_k}{f_k(z)-a_k},
\qquad
a_k=f_k(z_k).
\end{equation}
The interpolation condition at \(z_k\) has disappeared from the problem for \(f_{k+1}\).  Indeed, the remaining interpolation values are transformed according to
$$
w_j^{(k+1)}
 =
 \frac{z_j-z_k}{w_j^{(k)}-w_k^{(k)}},
 \qquad j>k,
$$
where \(w_j^{(0)}=w_j\), and
$$
a_k=w_k^{(k)}.
$$
These are precisely Thiele's reciprocal differences.

The crucial observation is that \eqref{thiele-reduction} can be inverted:
\begin{equation}\label{thiele-inverse}
f_k(z)
=a_k+\frac{z-z_k}{f_{k+1}(z)},
\end{equation}
which leads to the following
\begin{equation}\label{thiele-finite}
f(z)=
a_0+
\cfrac{z-z_0}{
a_1+
\cfrac{z-z_1}{
a_2+
\cfrac{z-z_2}{
\ddots+
\cfrac{z-z_{N-1}}{a_N+f_{N}(z)}}}},
\end{equation}
In this form, Thiele's procedure is naturally interpreted as a step-by-step reduction of an interpolation problem: at each stage one interpolation condition is removed, the remaining data are transformed, and the process is repeated. 
In particular, if we start with $N$ interpolation nodes \(z_0,z_1,\ldots,z_{N-1}\), then the function \(f_N\) is free of interpolation conditions and may therefore be chosen arbitrarily, yielding a parametrization of all solutions.

Even though Thiele's procedure produces a continued fraction, this fraction is not the main player for us. To introduce the object that will play that role, we consider what is known as the even part of a continued fraction. A standard reference for this construction is \cite{JT80}. We shall instead follow a matrix approach which, to the best of our knowledge, does not appear explicitly in the literature.
To this end, let us associate with the successive linear-fractional transformations in \eqref{thiele} the matrices
$$
M_j(z)
=
\begin{pmatrix}
0&z-z_{j-1}\\
1&a_j
\end{pmatrix},
\qquad j\geq 1.
$$
As usual, matrices differing by a nonzero scalar represent the same linear-fractional transformation. Therefore, \eqref{thiele} may be viewed formally as the infinite matrix product
$$
M_1(z)M_2(z)\cdots.
$$
The even part of \eqref{thiele} is the continued fraction whose successive convergents are the \(0\)-th, \(2\)-nd, \(4\)-th, \(\ldots\) convergents of the Thiele fraction. It is therefore natural to try to obtain it from the paired matrix product
$$
\bigl(M_1(z)M_2(z)\bigr)\bigl(M_3(z)M_4(z)\bigr)\cdots.
$$
There is, however, an immediate obstruction. Namely, the matrices $M_{2k-1}(z)M_{2k}(z)$ do not represent linear-fractional transformations of continued-fraction form, since their \((1,1)\)-entries need not vanish. Nevertheless, the idea is essentially correct and requires only a slight modification, as the following statement shows.

\begin{proposition}
Assume that \(a_{2k}\neq0\) for $k>0$. Then 
$$
M_{2k-1}(z)M_{2k}(z)=C_k(z)L_k(z),
$$
where
$$
C_k(z)
=
\begin{pmatrix}
0&a_{2k}(z-z_{2k-2})\\[1mm]
1&z-z_{2k-1}+a_{2k-1}a_{2k}
\end{pmatrix}
$$
and
$$
L_k(z)
=
\begin{pmatrix}
-\dfrac{z-z_{2k-1}}{a_{2k}}&0\\[2mm]
\dfrac1{a_{2k}}&1
\end{pmatrix}.
$$
As a result, we have that
$$
L_k(z)C_{k+1}(z)
=
\begin{pmatrix}
0&
-\dfrac{a_{2k+2}}{a_{2k}}
 (z-z_{2k-1})(z-z_{2k})
\\[3mm]
1&
z-z_{2k+1}
+a_{2k+1}a_{2k+2}
+\dfrac{a_{2k+2}}{a_{2k}}(z-z_{2k})
\end{pmatrix}.
$$
In particular, \(L_kC_{k+1}\) is again a continued-fraction matrix.
\end{proposition}

\begin{proof}
The factorization follows by direct multiplication:
$$
\begin{pmatrix}
0&z-z_{2k-2}\\
1&a_{2k-1}
\end{pmatrix}
\begin{pmatrix}
0&z-z_{2k-1}\\
1&a_{2k}
\end{pmatrix}
=
C_k(z)L_k(z).
$$
\end{proof}
The result of this modification is that
$$
M_1M_2\cdots M_{2n}\cdots 
=
C_1(L_1C_2)(L_2C_3)\cdots(L_{n-1}C_n)\cdots.
$$
The latter product generates the continued fraction of primary interest here, whose successive linear-fractional transformations are represented by the continued-fraction matrices $L_kC_{k+1}$. 
\begin{remark}
The matrix argument is not specific to the Thiele fraction. For a general continued fraction with step matrices
$$
M_j=
\begin{pmatrix}
0&a_j\\
1&b_j
\end{pmatrix},
$$
one has
$$
M_{2k-1}M_{2k}
=
\begin{pmatrix}
0&a_{2k-1}b_{2k}\\
1&a_{2k}+b_{2k-1}b_{2k}
\end{pmatrix}
\begin{pmatrix}
-\dfrac{a_{2k}}{b_{2k}}&0\\[2mm]
\dfrac1{b_{2k}}&1
\end{pmatrix}.
$$
Thus the same factorization, together with the elementary fact that the product of a lower triangular matrix and a continued-fraction matrix is again a continued-fraction matrix, gives the even contraction of an arbitrary generalized continued fraction. The factorization used here may be viewed as a scalar \(2\times2\) version of the triangular elimination underlying the Schur-complement construction. It is curious that the same triangular mechanism already present in Schur's work also provides a particularly simple way of obtaining the continued-fraction contraction.
\end{remark}

In order to write the continued fraction explicitly, let us set
$$
\rho_{k-1}=\frac{a_{2k+2}}{a_{2k}},\qquad k\geq1,
$$
and
$$
\beta_0z-\alpha_0=z-z_1+a_1a_2,
$$

$$
\beta_{k-1}z-\alpha_{k-1}
=
z-z_{2k-1}
+a_{2k-1}a_{2k}
+\kappa_k(z-z_{2k-2}),
\qquad k\geq2.
$$
Then the matrix product gives
\begin{equation}\label{eq:thiele-even-part}
\mathcal T_{\rm even}(z)
=
c_0+
\cfrac{c_2(z-z_0)}
{\beta_0z-\alpha_0
-
\cfrac{\rho_0(z-z_1)(z-z_2)}
{\beta_1z-\alpha_1
-
\cfrac{\rho_1(z-z_3)(z-z_4)}
{\beta_2z-\alpha_2-\ddots}}}.
\end{equation}
Finally, the change of variable \(z\mapsto 1/z\) preserves the structure of the contracted fraction \eqref{eq:thiele-even-part}, up to a normalization of its coefficients, provided none of the interpolation nodes is zero. If the first node is zero, the corresponding linear factor becomes constant and hence disappears after transforming back; this is precisely the case corresponding to the interpolation condition \(f(\infty)=0\). For simplicity, in what follows we assume that the first interpolation is always \(f(\infty)=0\) and so the fraction takes the form
\begin{equation}\label{eq:thiele-RII}
 \Phi(z)
 =
 \cfrac{1}{
 \beta_0z-\alpha_0
 -
 \cfrac{\rho_0(z-z_1)(z-z_2)}
 {\beta_1z-\alpha_1
 -
 \cfrac{\rho_1(z-z_3)(z-z_4)}
 {\beta_2z-\alpha_2-\ddots}}}.
\end{equation}
The continued fraction \eqref{eq:thiele-RII} is commonly referred to as an \(R_{II}\)-fraction. This terminology was introduced in \cite{IsmailMasson95}, although continued fractions of this form had appeared much earlier in various guises. The \(R_{II}\)-fraction, together with the recurrence relation it generates, will be the main object of interest in what follows. In particular, the denominator polynomials of its convergents satisfy
\begin{equation}\label{eq:even-RII}
P_{j+1}(z)
=(\beta_j z-\alpha_j)P_j(z)
-\rho_{j-1}
(z-z_{2j-1})(z-z_{2j})
P_{j-1}(z),
\qquad j\ge1,
\end{equation}
which is a recurrence relation of \(R_{II}\) type.

\section{Detour to Euler's formula}

There is another classical route to continued fractions which is useful for our purposes. It goes back to Euler and transforms a series into a continued fraction. The point of the present detour is purely algebraic. Namely, we will see that whenever the ratios of consecutive terms of the series are linear-fractional functions of the spectral variable, Euler's formula naturally produces an \(R_{II}\)-fraction. 

We begin by recalling Euler's formula. To this end, let
$$
S_n=s_0+s_1+\cdots+s_n
$$
and set
$$
r_k=\frac{s_k}{s_{k-1}},\qquad k=1,\ldots,n.
$$
Then
\begin{equation}\label{eq:euler-formula}
S_n
=
\cfrac{s_0}
{1-\cfrac{r_1}
{1+r_1-\cfrac{r_2}
{1+r_2-\ddots
-\cfrac{r_n}{1+r_n}}}}.
\end{equation}
For finite sums this is simply an algebraic identity. Its infinite counterpart is the familiar Euler transformation of a series into a continued fraction (see, for example, \cite[Section 4.1]{Khrushchev08}).

The relevance of \eqref{eq:euler-formula} to the present discussion becomes apparent already for a rational Nevanlinna function, which a rational function that has the admits the following representation:
\begin{equation}\label{eq:finite-discrete-m}
m_N(z)
=
\sum_{k=0}^{N}\frac{\mu_k}{t_k-z},
\qquad
\mu_k>0,\quad t_k\in\mathbb R,
\end{equation}
where the \(t_k\)'s are distinct. Setting
\begin{equation}\label{eq:euler-ratio}
\gamma_k=\frac{\mu_k}{\mu_{k-1}},
\qquad k=1,\ldots,N,
\end{equation}
for the summands in \eqref{eq:finite-discrete-m} we get
$$
r_k(z)
=
\gamma_k\frac{z-t_{k-1}}{z-t_k}.
$$
Substitution into Euler's formula, followed by elementary equivalence transformations, gives
\begin{equation}\label{eq:euler-RII-finite}
m_N(z)
=-\cfrac{\mu_0}
{z-t_0
-\cfrac{\gamma_1(z-t_0)^2}
{d_1(z)
-\cfrac{\gamma_2(z-t_1)^2}
{d_2(z)
-\ddots
-\cfrac{\gamma_N(z-t_{N-1})^2}
{d_N(z)}}}},
\end{equation}
where
\begin{equation}\label{eq:euler-dk}
d_k(z)
=
(1+\gamma_k)z-
\bigl(t_k+\gamma_k t_{k-1}\bigr).
\end{equation}
Thus Euler's formula transforms a partial-fraction expansion into a finite \(R_{II}\)-fraction. It is a degenerate, or confluent, \(R_{II}\)-fraction in the sense that the two linear factors in every quadratic numerator coincide (cf. formula \eqref{eq:thiele-RII}).

\begin{remark}
It is worth noting that exactly the same reduction appears in \cite{AglerYoung} in the context of boundary interpolation. If a Nevanlinna function \(f\) has a pole at \(x\in\mathbb R\) with residue \(R\), its reduction at \(x\) is defined by
\[
g(\lambda)
=
f(z)-\frac{R}{z-x}.
\]
For the rational Nevanlinna functions considered above,
\[
m_N(\lambda)
=
\sum_{k=0}^N\frac{\mu_k}{t_k-z},
\]
the residue at \(t_N\) is \(-\mu_N\), and hence reduction at \(t_N\) gives precisely
\[
m_N(z)\longmapsto m_{N-1}(z).
\]
Thus, at a pole, the reduction from \cite{AglerYoung} is exactly the elementary pole-removal step underlying Euler's construction.

In \cite{AglerYoung} the same reduction procedure is also extended to points at which the function is analytic, and is used as a tool for boundary Nevanlinna--Pick interpolation. The emphasis there is not on producing a continued fraction, although repeated reduction and augmentation are naturally encoded by linear-fractional transformations and their \(2\times2\) matrix products.
\end{remark}

As is known, a continued fraction generates a recurrence relation (for example, see \cite[Section 2.1]{JT80}) and in our case the corresponding recurrence has an especially simple structure. Namely, letting
$$
P_0(z)=1,\qquad P_1(z)=z-t_0,
$$
the standard formula for the denominators of convergents to a continued fractions applied to \eqref{eq:euler-RII-finite} yields 
\begin{equation}\label{eq:euler-RII-rec}
P_{k+1}(z)=d_k(z)P_k(z)-
\gamma_k(z-t_{k-1})^2P_{k-1}(z),
\qquad k=1,\dots,N.
\end{equation}
The three-term recurrence relation \eqref{eq:euler-RII-rec} can be symmetrized as follows

As is known, a continued fraction generates a recurrence relation (for example, see \cite[Section~2.1]{JT80}), and in our case the corresponding recurrence has an especially simple structure. Namely, letting
$$
P_0(z)=1,\qquad P_1(z)=z-t_0,
$$
The three-term recurrence relation \eqref{eq:euler-RII-rec} can be symmetrized in the same way as the \(R_{II}\)-recurrences considered in \cite{DZ09}. To this end, introduce
$$
\widehat P_0(z)=1,\qquad \widehat P_k(z)
=
\frac{P_k(z)}
{\sqrt{\gamma_1\cdots\gamma_k}\prod_{j=0}^{k-1}(t_j-z)},
\qquad k=1,...,N+1,
$$
where we set $\gamma_{N+1}=1$ for convenience.
Then \eqref{eq:euler-RII-rec} takes the symmetric form
\begin{equation}\label{eq:euler-symmetric-rec}
\sqrt{\gamma_{k+1}}(z-t_k)\widehat P_{k+1}(z)
+
((1+\gamma_k)z-
\bigl(t_k+\gamma_k t_{k-1}\bigr))
\widehat P_k(z)
+
\sqrt{\gamma_k}(z-t_{k-1})\widehat P_{k-1}(z)
=0.
\end{equation}
It is important to emphasize that the normalization above is applied to the recurrence relation itself rather than to a particular solution. Indeed, the recurrence is the object that encodes the continued fraction and, after symmetrization, gives rise naturally to the corresponding generalized eigenvalue problem. In the present special case, the distinguished denominator solution
\[
P_k(z)=\prod_{j=0}^{k-1}(z-t_j)
\]
is completely factorized, and therefore the transformation
\[
\widehat P_k(z)
=
\frac{P_k(z)}
{\kappa_k\prod_{j=0}^{k-1}(t_j-z)}
\]
reduces this particular solution to the constant sequence
\[
\widehat P_k(z)=\frac{(-1)^k}{\kappa_k}.
\]
This degeneracy reflects the especially simple structure of the Euler \(R_{II}\)-fraction. The role of the transformation is nevertheless to symmetrize the recurrence and reveal the associated pencil.
Namely, setting $\pi(z)=(\widehat P_0(z),\dots, \widehat P_{N}(z))^{\top}$ the recurrence \eqref{eq:euler-symmetric-rec} can be written as the generalized eigenvalue problem 
$$
J_N^{(1)}\pi(z)=zJ_N^{(2)}\pi(z),
$$
where
\begin{equation}\label{eq:rational-pencil-j2}
J_N^{(2)}
=
\begin{pmatrix}
1 & \sqrt{\gamma_1} & & & \\
\sqrt{\gamma_1} & 1+\gamma_1 & \sqrt{\gamma_2} & & \\
& \sqrt{\gamma_2} & 1+\gamma_2 & \ddots & \\
& & \ddots & \ddots & \sqrt{\gamma_N}\\
& & & \sqrt{\gamma_N} & 1+\gamma_N
\end{pmatrix}
\end{equation}
and
\begin{equation}\label{eq:rational-pencil-j1}
J_N^{(1)}
=
\begin{pmatrix}
t_0 & \sqrt{\gamma_1}t_0 & & & \\
\sqrt{\gamma_1}t_0
&
t_1+\gamma_1t_0
&
\sqrt{\gamma_2}t_1
&
&\\
&
\sqrt{\gamma_2}t_1
&
t_2+\gamma_2t_1
&
\ddots
&\\
&
&
\ddots
&
\ddots
&
\sqrt{\gamma_N}t_{N-1}\\
&
&
&
\sqrt{\gamma_N}t_{N-1}
&
t_N+\gamma_Nt_{N-1}
\end{pmatrix}.
\end{equation}
As in the case considered in \cite{DZ09}, the matrix $J_N^{(2)}$ has the Cholesky factorization. Indeed, if
$$
L_N=
\begin{pmatrix}
1 & & & &\\
\sqrt{\gamma_1}&1&&&\\
&\sqrt{\gamma_2}&1&&\\
&&\ddots&\ddots&\\
&&&\sqrt{\gamma_N}&1
\end{pmatrix}
$$
then
\begin{equation*}
J_N^{(2)}=L_NL_N^*.
\end{equation*}
In addition, in the current settings, we have that
\[
J_N^{(1)}=L_ND_NL_N^*, \qquad
D_N=\operatorname{diag}(t_0,t_1,\ldots,t_N).
\]
Hence,
\begin{equation}\label{eq:euler-pencil-factor}
J_N^{(1)}-zJ_N^{(2)}
=
L_N(D_N-zI)L_N^*.
\end{equation}
In particular, \(J_N^{(2)}\) is positive definite and the generalized eigenvalues of the pencil \(J_N^{(1)}-zJ_N^{(2)}\) are precisely the poles of the function $m_N$. 

The factorization \eqref{eq:euler-pencil-factor} also gives an immediate realization of the rational Nevanlinna function from which the construction started.

\begin{proposition}\label{prop:euler-realization}
Let
$$
m_N(z)
=
\sum_{j=0}^{N}\frac{\mu_j}{t_j-z},
\qquad
\mu_j>0.
$$
Then the matrices \(J_N^{(1)}\) and \(J_N^{(2)}\) defined by \eqref{eq:rational-pencil-j2}, \eqref{eq:rational-pencil-j1} and \eqref{eq:euler-ratio} satisfy
\begin{equation}\label{eq:euler-realization}
m_N(z)
=
\mu_0
\left(
\left(J_N^{(1)}-zJ_N^{(2)}\right)^{-1}e_0,e_0
\right)_{\mathbb{C}^{N+1}},
\end{equation}
where $e_0=(1,0,\dots,0)^{\top}$ is the first standard basis vector. 
Thus the \(R_{II}\)-fraction obtained from Euler's formula yields a finite-dimensional realization of \(m_N\) by a pencil of Jacobi matrices. In this realization, \(m_N\) will be referred to as the \(m\)-function of the pencil.
\end{proposition}

\begin{proof}
For $z\in\mathbb{C}\setminus\{t_0,t_1,\dots,t_N\}$, from \eqref{eq:euler-pencil-factor} we see that
$$
\left(J_N^{(1)}-zJ_N^{(2)}\right)^{-1}
=
L_N^{-*}(D_N-zI)^{-1}L_N^{-1}
$$
since $L_N$ is invertible.
Moreover, $L_N$ is bidiagonal and one can find its inverse, which in particular yields
$$
L_N^{-1}e_0
=
\begin{pmatrix}
1\\[1mm]
-\sqrt{\gamma_1}\\
\sqrt{\gamma_1\gamma_2}\\
\vdots\\
(-1)^N\sqrt{\gamma_1\cdots\gamma_N}
\end{pmatrix}.
$$
Next, we have that $\gamma_1\cdots\gamma_j=\dfrac{\mu_j}{\mu_0}$ and, therefore, we get
$$
\mu_0
\left(
\left(J_N^{(1)}-zJ_N^{(2)}\right)^{-1}e_0,e_0
\right)
=
\sum_{j=0}^{N}\frac{\mu_j}{t_j-z}
=
m_N(z).
$$
\end{proof}
\begin{remark}
The realization \eqref{eq:euler-realization} can also be understood from a more elementary point of view. Indeed, even without introducing the tridiagonal pencil, one can write the rational function \(m_N\) as
\[
m_N(z)
=
v_N^{\top}(D_N-zI)^{-1}v_N,
\]
where
\[
D_N=\operatorname{diag}(t_0,\ldots,t_N),
\qquad
v_N=
\left(
\sqrt{\mu_0},
-\sqrt{\mu_1},
\sqrt{\mu_2},
\ldots,
(-1)^N\sqrt{\mu_N}
\right)^{\top}.
\]
This is simply the pole-residue representation written in matrix form. Thus the finite-dimensional realization itself is immediate. The point of the Euler construction is that it transforms this diagonal realization into the structured tridiagonal pencil associated with the \(R_{II}\)-fraction. 
\end{remark}

The same construction can be applied to some meromorphic functions. In particular, let us consider one of the classical weights, the Meixner weight
$$
\mu_k=\frac{(1-\alpha)_k}{k!}q^k,
\qquad
0<q<1,\qquad \alpha<1,
$$
where 
\[
(a)_k=a(a+1)\cdots(a+k-1),\qquad k\geq1,
\]
denotes the Pochhammer symbol (or rising factorial), with
\[
(a)_0=1.
\]
The corresponding meromorphic Nevanlinna function has the form
\begin{equation}\label{eq:meixner-m}
m_{\alpha,q}(z)
=\sum_{k=0}^{\infty}
\frac{(1-\alpha)_kq^k}{k!(k-z)}.
\end{equation}
For \eqref{eq:meixner-m}, the ratio of two consecutive summands is
$$
r_k(z)
=
q\,\frac{k-\alpha}{k}\,
\frac{z-k+1}{z-k}.
$$
Thus, with
$$
\gamma_k
=
q\frac{k-\alpha}{k},
$$
Euler's formula gives formally
\begin{equation}\label{eq:meixner-RII}
m_{\alpha,q}(z)
=-\cfrac{1}
{z
-\cfrac{\gamma_1z^2}
{d_1(z)
-\cfrac{\gamma_2(z-1)^2}
{d_2(z)
-\cfrac{\gamma_3(z-2)^2}
{d_3(z)-\ddots}}}},
\end{equation}
where
$$
d_k(z)
=
(1+\gamma_k)z
-
\bigl(k+\gamma_k(k-1)\bigr).
$$
Thus the Meixner case leads to an infinite \(R_{II}\)-fraction of exactly the same degenerate type. At the purely formal level, the associated recurrence is
\begin{equation}\label{eq:meixner-RII-rec}
P_{k+1}(z)
=d_k(z)P_k(z)
-\gamma_k(z-k+1)^2P_{k-1}(z),
\end{equation}
and its denominator polynomial solution is simply
$$
P_{k+1}(z)=z(z-1)\cdots(z-k).
$$
Similarly, the corresponding pencil \(J^{(1)}-zJ^{(2)}\) is determined by
\[
J^{(2)}
=
\begin{pmatrix}
1
&
\sqrt{q(1-\alpha)}
&
&
&
\\[1mm]
\sqrt{q(1-\alpha)}
&
1+q(1-\alpha)
&
\sqrt{\dfrac{q(2-\alpha)}{2}}
&
&
\\[3mm]
&
\sqrt{\dfrac{q(2-\alpha)}{2}}
&
1+\dfrac{q(2-\alpha)}{2}
&
\sqrt{\dfrac{q(3-\alpha)}{3}}
&
\\
&
&
\ddots&\ddots&\ddots
\end{pmatrix}
\]
and
\[
J^{(1)}
=
\begin{pmatrix}
0&0&&&\\[1mm]
0&1&
\sqrt{\dfrac{q(2-\alpha)}{2}}
&&\\[3mm]
&
\sqrt{\dfrac{q(2-\alpha)}{2}}
&
2+\dfrac{q(2-\alpha)}{2}
&
2\sqrt{\dfrac{q(3-\alpha)}{3}}
&
\\[3mm]
&
&
2\sqrt{\dfrac{q(3-\alpha)}{3}}
&
3+\dfrac{2q(3-\alpha)}{3}
&
\ddots
\\
&
&
&
\ddots&\ddots
\end{pmatrix}.
\]
As in the finite case, if
\[
L=
\begin{pmatrix}
1&&&&\\
\sqrt{q(1-\alpha)}&1&&&\\
&\sqrt{\dfrac{q(2-\alpha)}{2}}&1&&\\
&&\sqrt{\dfrac{q(3-\alpha)}{3}}&1&\\
&&&\ddots&\ddots
\end{pmatrix},
\qquad
D=\operatorname{diag}(0,1,2,\ldots),
\]
then
\[
J^{(2)}=LL^*,
\qquad
J^{(1)}=LDL^*,
\]
and therefore
\[
J^{(1)}-zJ^{(2)}
=
L(D-zI)L^*.
\]
Since $L=I+A$ with $\Vert A\Vert<1$ we conclude that $L$ is invertible in $\ell^2$ and, hence,
\[
m_{\alpha,q}(z)
=
\mu_0
\left(
\left(J^{(1)}-zJ^{(2)}\right)^{-1}e_0,e_0
\right)_{\ell^2}
\]
for $z\in\mathbb{C}\setminus\mathbb{Z}_+$.
\begin{remark}
The observation explored in this section is not specific to Nevanlinna functions and applies, for example, to a broad class of hypergeometric functions. Recall that the generalized hypergeometric function is defined by
\[
{}_pF_q\!\left(
\begin{matrix}
a_1,\ldots,a_p\\
b_1,\ldots,b_q
\end{matrix};\xi
\right)
=
\sum_{k=0}^{\infty}
\frac{(a_1)_k\cdots(a_p)_k}
     {(b_1)_k\cdots(b_q)_k}
\frac{\xi^k}{k!},
\]
whenever the series is defined. Consider now the particular case
\[
{}_pF_q\!\left(
\begin{matrix}
A-z,a_2,\ldots,a_p\\
B-z,b_2,\ldots,b_q
\end{matrix};\xi
\right)
=
\sum_{k\ge0}T_k(z).
\]
Then
\[
\frac{T_{k+1}(z)}{T_k(z)}
=
c_k\,\frac{z-A-k}{z-B-k},
\]
where
\[
c_k
=
\frac{\xi}{k+1}
\frac{(a_2+k)\cdots(a_p+k)}
     {(b_2+k)\cdots(b_q+k)}.
\]
Thus the ratios of consecutive terms are linear-fractional functions of the spectral variable. Applying Euler's formula and performing the same equivalence transformations as above, we obtain formally
\[
{}_pF_q\!\left(
\begin{matrix}
A-z,a_2,\ldots,a_p\\
B-z,b_2,\ldots,b_q
\end{matrix};\xi
\right)
=
\cfrac{1}
{1-
 \cfrac{c_0(z-A)}
 {d_0(z)-
  \cfrac{c_1(z-B)(z-A-1)}
  {d_1(z)-
   \cfrac{c_2(z-B-1)(z-A-2)}
   {d_2(z)-\ddots}}}},
\]
where
\[
d_k(z)
=
(z-B-k)+c_k(z-A-k).
\]
Hence Euler's transformation produces, at least formally, a continued fraction of \(R_{II}\) type. In particular, the appearance of \(R_{II}\)-fractions in hypergeometric examples is a consequence of the linear-fractional dependence of the ratios of consecutive terms on the spectral variable, rather than a feature of any particular hypergeometric family.
\end{remark}

\section{From the Schur algorithm to Wall pencils}

As problems arising in analysis and its applications increasingly involved analytic functions subject to additional constraints, considerable attention turned to distinguished classes of analytic functions. Moment problems, coefficient problems, and interpolation problems naturally led to settings in which boundedness, positivity, or related analytic properties had to be preserved. In particular, the study of analytic functions bounded in the unit disk led Schur, in 1917, to introduce a recursive procedure for describing functions \(f\) analytic in the unit disk and satisfying \(|f(z)|\leq1\). The class of such functions is now known as the Schur class, and the corresponding recursive procedure as the Schur algorithm. Its essential feature is that the problem is reduced step by step while remaining within the same class of functions.

Unlike the constructions considered in the preceding sections, the original Schur algorithm is not, strictly speaking, a continued-fraction algorithm. Its defining feature is a sequence of linear-fractional transformations preserving the Schur class. Nevertheless, precisely because these transformations are linear-fractional, the algorithm can be reorganized into a continued fraction. Historically, the passage from Schur's recursion to a genuine continued-fraction representation took several additional steps.

First, recall that if \(f=f_0\) belongs to the Schur class, then the Schur algorithm is defined recursively by
\begin{equation}\label{eq:schur-alg}
 f_{n+1}(z)
 =
 \frac{f_n(z)-\gamma_n}
 {z\bigl(1-\overline{\gamma_n}f_n(z)\bigr)},
 \qquad
 \gamma_n=f_n(0),\qquad |\gamma_n|<1.
\end{equation}
Equivalently,
\begin{equation}\label{eq:schur-backward}
 f_n(z)
 =
 \frac{\gamma_n+zf_{n+1}(z)}
 {1+\overline{\gamma_n}zf_{n+1}(z)}.
\end{equation}
Schur himself referred to this procedure as ``continued fraction-like.'' Shortly after Schur's work, Hamel \cite{Hamel}, apparently independently, developed a closely related recursive construction for analytic functions bounded in the unit disk and expressed it in continued-fraction form. 

Wall subsequently returned to Schur's algorithm and transformed it into an actual continued-fraction algorithm. The even part of Wall's continued fraction gives a simpler contracted fraction, referred to by Khrushchev as the Geronimus continued fraction. Geronimus obtained this representation independently at approximately the same time as Wall, while a continued fraction essentially of the same form had already appeared in the earlier work of Hamel. After passing from the unit disk to the upper half-plane by the Cayley transform, this contracted fraction takes the form of an \(R_{II}\)-fraction, with quadratic partial numerators proportional to
\[
\lambda^2+1=(\lambda-i)(\lambda+i).
\]
These constructions are discussed in Wall's monograph \cite{Wall48}; see also \cite{Khrushchev08} for a modern account.

We now describe this development using matrix representations of linear-fractional transformations and continued fractions. First, notice that the backward Schur transformation \eqref{eq:schur-backward} is represented by
\[
S_n(z)
=
\begin{pmatrix}
z&\gamma_n\\
\overline{\gamma_n}z&1
\end{pmatrix}.
\]
Equivalently,
\[
f_n(z)
=
\gamma_n+
\cfrac{(1-|\gamma_n|^2)z}
{\overline{\gamma_n}z+\cfrac{1}{f_{n+1}(z)}}.
\]
Iterating this identity gives Wall's continued fraction
\begin{equation}\label{eq:wall-first}
f(z)
=
\gamma_0+
\cfrac{(1-|\gamma_0|^2)z}
{\overline{\gamma_0}z+
 \cfrac{1}
 {\gamma_1+
  \cfrac{(1-|\gamma_1|^2)z}
  {\overline{\gamma_1}z+
   \cfrac{1}
   {\gamma_2+\ddots}}}}.
\end{equation}
This is the original, uncontracted Wall fraction; see \cite{Wall48,Khrushchev08}.

Let us now recover its contraction directly from the matrices. Assume for simplicity that \(\gamma_n\neq0\). We have
\[
S_n(z)=C_nL_n(z),
\]
where
\[
C_n=
\begin{pmatrix}
0&\gamma_n\\
1&1
\end{pmatrix},
\qquad
L_n(z)=
\begin{pmatrix}
-\dfrac{1-|\gamma_n|^2}{\gamma_n}z&0\\[2mm]
\dfrac{z}{\gamma_n}&1
\end{pmatrix}.
\]
Consequently,
\[
S_0S_1\cdots S_N
=
C_0L_0C_1L_1\cdots C_NL_N
=
C_0(L_0C_1)(L_1C_2)\cdots(L_{N-1}C_N)L_N.
\]
The first product corresponds to the successive Schur steps underlying the Wall fraction \eqref{eq:wall-first}. Regrouping the matrices in the second form amounts precisely to contracting every two successive levels of that continued fraction, and therefore produces its even part.
Indeed, introducing
\[
\rho_n=\frac{\gamma_{n+1}}{\gamma_n},
\qquad
\sigma_n=1-|\gamma_n|^2,
\]
we obtain
\[
L_n(z)C_{n+1}
=
\begin{pmatrix}
0&-\sigma_n\rho_n z\\
1&1+\rho_nz
\end{pmatrix}.
\]
Thus, apart from the terminal factor, the regrouped product consists entirely of continued-fraction matrices and gives
\begin{equation}\label{eq:hamel-geronimus}
f(z)
=
\cfrac{\gamma_0}
{1-
 \cfrac{\sigma_0\rho_0 z}
 {1+\rho_0z-
  \cfrac{\sigma_1\rho_1 z}
  {1+\rho_1z-
   \cfrac{\sigma_2\rho_2 z}
   {1+\rho_2z-\ddots}}}}.
\end{equation}
This is the continued fraction usually associated with Geronimus and, up to inessential differences of normalization, with the earlier construction of Hamel. In the terminology of \cite{Khrushchev08}, it is the even part of Wall's continued fraction \eqref{eq:wall-first}. Thus the passage from the Wall fraction to the continued fraction considered by Hamel and Geronimus appears, at the matrix level, simply as the regrouping
\[
C_0L_0C_1L_1\cdots
\quad\longmapsto\quad
C_0(L_0C_1)(L_1C_2)\cdots .
\]

If we now apply the Cayley transformation
\[
z=\frac{i-\lambda}{i+\lambda}
\]
to \eqref{eq:hamel-geronimus} and perform the usual equivalence transformations of continued fractions, we obtain
\begin{equation}\label{eq:hamel-geronimus-RII}
f\left(\frac{i-\lambda}{i+\lambda}\right)
=
\cfrac{\gamma_0}{
1-
\cfrac{\sigma_0\rho_0(i-\lambda)}
{\beta_0(\lambda)+
 \cfrac{\sigma_1\rho_1(\lambda^2+1)}
 {\beta_1(\lambda)+
  \cfrac{\sigma_2\rho_2(\lambda^2+1)}
  {\beta_2(\lambda)+\ddots}}}},
\end{equation}
where
\[
\beta_n(\lambda)
=
(1-\rho_n)\lambda+i(1+\rho_n).
\]
Thus, after the initial linear numerator, all partial numerators contain the fixed quadratic factor
\[
\lambda^2+1=(\lambda-i)(\lambda+i).
\]
In particular, \eqref{eq:hamel-geronimus-RII} already has the characteristic \(R_{II}\) structure.

The same transformation is especially transparent at the matrix level. Under the Cayley transformation, and using the fact that multiplication of a matrix representing a linear-fractional transformation by a nonzero scalar does not change the transformation, the continued-fraction matrix
\[
L_n(z)C_{n+1}
=
\begin{pmatrix}
0&-\sigma_n\rho_nz\\
1&1+\rho_nz
\end{pmatrix}
\]
may be replaced by
\[
\widetilde W_n(\lambda)
=
\begin{pmatrix}
0&-\sigma_n\rho_n(i-\lambda)\\
i+\lambda&(i+\lambda)+\rho_n(i-\lambda)
\end{pmatrix}.
\]
Write
\[
\widetilde W_n(\lambda)
=
R_n(\lambda)D(\lambda),
\]
where
\[
R_n(\lambda)
=
\begin{pmatrix}
0&-\sigma_n\rho_n(i-\lambda)\\
1&\beta_n(\lambda)
\end{pmatrix},
\qquad
D(\lambda)
=
\begin{pmatrix}
i+\lambda&0\\
0&1
\end{pmatrix},
\]
and
\[
\beta_n(\lambda)
=
(i+\lambda)+\rho_n(i-\lambda)
=
(1-\rho_n)\lambda+i(1+\rho_n).
\]
The diagonal factor can then be absorbed into the next continued-fraction matrix:
\[
D(\lambda)R_{n+1}(\lambda)
=
\begin{pmatrix}
0&
-\sigma_{n+1}\rho_{n+1}(i+\lambda)(i-\lambda)
\\[1mm]
1&\beta_{n+1}(\lambda)
\end{pmatrix}.
\]
Since
\[
-(i+\lambda)(i-\lambda)=\lambda^2+1,
\]
this becomes
\[
D(\lambda)R_{n+1}(\lambda)
=
\begin{pmatrix}
0&
\sigma_{n+1}\rho_{n+1}(\lambda^2+1)
\\
1&
(1-\rho_{n+1})\lambda+i(1+\rho_{n+1})
\end{pmatrix}.
\]
Consequently,
\[
\widetilde W_0\widetilde W_1\cdots\widetilde W_N
=
R_0
(DR_1)(DR_2)\cdots(DR_N)D.
\]
Thus, after the initial factor \(R_0\), every factor in the product is a continued-fraction matrix of \(R_{II}\) type, with the same pair of nodes \(i\) and \(-i\). The terminal diagonal factor \(D\) fixes \(0\) and therefore does not affect the corresponding convergents.

Although \eqref{eq:hamel-geronimus-RII} is a perfectly legitimate continued fraction, it does not seem to have played an independent role in the theory. One drawback is that its individual linear-fractional transformations do not preserve a distinguished analytic class. Wall's further normalization remedies this difficulty and leads to a continued fraction with real parameters.

There is a convenient matrix interpretation of this normalization. Suppose that the successive steps of the Schur algorithm are represented by
\[
S_n(z)
=
\begin{pmatrix}
z&\gamma_n\\
\overline{\gamma_n}z&1
\end{pmatrix}.
\]
A change of the intermediate linear-fractional coordinates amounts to replacing \(S_n\) by
\[
\widetilde S_n(z)
=
G_n(z)^{-1}S_n(z)G_{n+1}(z),
\]
where \(G_n(z)\) is invertible. Indeed,
\[
\widetilde S_0\widetilde S_1\cdots\widetilde S_N
=
G_0^{-1}S_0S_1\cdots S_NG_{N+1},
\]
so all intermediate changes of coordinates cancel.

Wall's choice corresponds to
\[
G_n(z)^{-1}
=
\begin{pmatrix}
-\delta_n&1\\
\delta_n z&1
\end{pmatrix},
\qquad |\delta_n|=1,
\]
or, equivalently,
\[
h_n(z)
=
\frac{1-\delta_n f_n(z)}
     {1+\delta_n zf_n(z)}.
\]
The phase \(\delta_{n+1}\) is chosen so that the \((1,1)\)-entry of
\[
\widetilde S_n(z)=G_n(z)^{-1}S_n(z)G_{n+1}(z)
\]
vanishes. This gives
\[
\delta_{n+1}
=
\frac{\overline{\gamma_n}-\delta_n}
     {1-\gamma_n\delta_n},
\qquad
\delta_0=1,
\]
and \(|\delta_{n+1}|=1\) whenever \(|\delta_n|=1\).

With this choice, and up to multiplication by a nonzero scalar,
\[
\widetilde S_n(z)
\sim
\begin{pmatrix}
0&
|1-\gamma_n\delta_n|^2
\\[1mm]
(1-|\gamma_n|^2)z&
(1-\overline{\gamma_n}\overline{\delta_n})
-(1-\gamma_n\delta_n)z
\end{pmatrix}.
\]
Put
\[
a_n=1-\gamma_n\delta_n
\]
and introduce Wall's real parameters
\[
g_{n+1}
=
\frac{|a_n|^2}{2\Re a_n},
\qquad
r_{n+1}
=
-\frac{\Im a_n}{\Re a_n}.
\]
Since
\[
1-|\gamma_n|^2
=
2\Re a_n-|a_n|^2,
\]
the transformed matrix takes the particularly simple form
\[
\widetilde S_n(z)
\sim
\begin{pmatrix}
0&2g_{n+1}\\[1mm]
2(1-g_{n+1})z&
(1+ir_{n+1})-(1-ir_{n+1})z
\end{pmatrix},
\]
where
\[
0<g_{n+1}<1,
\qquad
r_{n+1}\in\mathbb R.
\]
Thus the parameters \(g_{n+1}\) and \(r_{n+1}\) may be viewed as real coordinates of a distinguished gauge representative of the \(n\)-th Schur step.

Since \(\delta_0=1\),
\[
h_0(z)=\frac{1-f(z)}{1+zf(z)},
\]
and iteration of the preceding transformations gives
\[
h_0(z)
=
\cfrac{2g_1}{
(1+ir_1)-(1-ir_1)z
+
\cfrac{4(1-g_1)g_2z}{
(1+ir_2)-(1-ir_2)z
+
\cfrac{4(1-g_2)g_3z}{
(1+ir_3)-(1-ir_3)z+\ddots}}}.
\]
This is the normalized continued fraction underlying Wall's subsequent passage to Carath\'eodory and Nevanlinna functions.

\begin{theorem}[\cite{D17},\cite{Wall48}]\label{thm:wall-cf}
Let \(f\) be a Schur function, and let \(g_k,r_k\) be the real parameters obtained from Wall's normalization of its Schur parameters. Define
\[
F(z)=\frac{1+zf(z)}{1-zf(z)}
\]
and
\[
\varphi(\lambda)
=
iF\!\left(\frac{i-\lambda}{i+\lambda}\right).
\]
Then
\begin{equation}\label{eq:wall-RII}
\varphi(\lambda)
=
-\cfrac{1}{
\lambda-
\cfrac{g_1(\lambda^2+1)}
{\lambda-r_1-
\cfrac{(1-g_1)g_2(\lambda^2+1)}
{\lambda-r_2-
\cfrac{(1-g_2)g_3(\lambda^2+1)}
{\lambda-r_3-\ddots}}}}.
\end{equation}
Here
\[
0<g_k<1,
\qquad
r_k\in\mathbb R.
\]
The equality is understood in the formal continued-fraction sense: every finite truncation of the right-hand side is obtained from the corresponding finite product of the continued-fraction matrices arising in Wall's construction. In particular, \eqref{eq:wall-RII} is an \(R_{II}\)-fraction with the fixed pair of nodes \(i\) and \(-i\).
\end{theorem}

As usual, the continued fraction \eqref{eq:wall-RII} generates a three-term recurrence. Setting
\[
g_0=r_0=0,
\]
one obtains
\begin{equation}\label{eq:wall-rec}
u_{j+1}(\lambda)
-(\lambda-r_j)u_j(\lambda)
+(1-g_{j-1})g_j(\lambda^2+1)u_{j-1}(\lambda)
=0.
\end{equation}
Although \eqref{eq:wall-rec} is quadratic in the spectral parameter, it admits a particularly useful renormalization. Put
\[
\mathfrak a_j=r_j,
\qquad
\mathfrak b_j=\sqrt{(1-g_j)g_{j+1}},
\]
and define
\[
\widehat u_0=u_0,
\qquad
\widehat u_j(\lambda)
=
\frac{u_j(\lambda)}
{\mathfrak b_0\cdots\mathfrak b_{j-1}(i-\lambda)^j},
\qquad j\ge1.
\]
Then \eqref{eq:wall-rec} becomes
\begin{equation}\label{eq:wall-symmetric-rec}
\mathfrak b_j(i-\lambda)\widehat u_{j+1}
-(\lambda-\mathfrak a_j)\widehat u_j
-\mathfrak b_{j-1}(i+\lambda)\widehat u_{j-1}
=0.
\end{equation}
Consequently, \eqref{eq:wall-symmetric-rec} is the generalized eigenvalue problem
\begin{equation}\label{eq:wall-pencil}
(H-\lambda J)\widehat u=0,
\end{equation}
where
\[
H=
\begin{pmatrix}
\mathfrak a_0&i\mathfrak b_0&&\\
-i\mathfrak b_0&\mathfrak a_1&i\mathfrak b_1&\\
&-i\mathfrak b_1&\mathfrak a_2&\ddots\\
&&\ddots&\ddots
\end{pmatrix},
\qquad
J=
\begin{pmatrix}
1&\mathfrak b_0&&\\
\mathfrak b_0&1&\mathfrak b_1&\\
&\mathfrak b_1&1&\ddots\\
&&\ddots&\ddots
\end{pmatrix}.
\]
Both matrices are Hermitian Jacobi matrices, while the special form
\[
\mathfrak b_j^2=(1-g_j)g_{j+1}
\]
shows that the off-diagonal coefficients of \(J\) form a chain sequence. In particular, \(J\) is nonnegative. Following \cite{D17}, we shall refer to the linear pencils \(H-\lambda J\) obtained in this way as \emph{Wall pencils}. Thus Wall's modification of the Schur algorithm leads, through the \(R_{II}\)-fraction \eqref{eq:wall-RII}, directly to a spectral problem for a pair of Jacobi matrices.

\section{A detour to CMV matrices}

We now take a short detour to CMV matrices. Perhaps somewhat surprisingly, the same general strategy that underlies Wall's construction can be reorganized to lead to the CMV setting as well. We shall show how this connection emerges from the Schur transformations.

As before, let
\[
f_n(z)
=
\frac{\alpha_n+zf_{n+1}(z)}
     {1+\overline{\alpha_n}zf_{n+1}(z)},
\qquad
\rho_n^2=1-|\alpha_n|^2,
\]
represent the backward Schur transformation by
\[
S_n(z)
=
\begin{pmatrix}
z&\alpha_n\\
\overline{\alpha_n}z&1
\end{pmatrix}.
\]
We first make the the transformation 
\[
f\longmapsto \frac{f}{1-zf}.
\]
which corresponds to
\[
H(z)=
\begin{pmatrix}
1&0\\
-z&1
\end{pmatrix}.
\]
A direct calculation gives
\[
H(z)S_0(z)
=
W_0(z)R_0(z),
\]
where
\[
W_0(z)
=
\begin{pmatrix}
0&1\\
1&\overline{\alpha_0}-z
\end{pmatrix},
\qquad
R_0(z)
=
\begin{pmatrix}
0&\rho_0^2\\
z&\alpha_0
\end{pmatrix}.
\]
The remaining factorizations alternate between two forms. For \(k\geq0\), put
\[
R_{2k}(z)
=
\begin{pmatrix}
0&\rho_{2k}^2\\
z&\alpha_{2k}
\end{pmatrix}.
\]
Then
\begin{equation}\label{eq:cmv-factor-odd}
R_{2k}(z)S_{2k+1}(z)
=
W_{2k+1}(z)R_{2k+1}(z),
\end{equation}
where
\[
W_{2k+1}(z)
=
\begin{pmatrix}
0&\rho_{2k}^2\\
1&\alpha_{2k}+\alpha_{2k+1}z
\end{pmatrix}
\]
and
\[
R_{2k+1}(z)
=
\begin{pmatrix}
\rho_{2k+1}^2z^2&0\\
\overline{\alpha_{2k+1}}z&1
\end{pmatrix}.
\]
Indeed,
\[
R_{2k}S_{2k+1}
=
\begin{pmatrix}
\rho_{2k}^2\overline{\alpha_{2k+1}}z&
\rho_{2k}^2\\[1mm]
z^2+\alpha_{2k}\overline{\alpha_{2k+1}}z&
\alpha_{2k}+\alpha_{2k+1}z
\end{pmatrix},
\]
and \eqref{eq:cmv-factor-odd} follows from
\[
\rho_{2k+1}^2
+
|\alpha_{2k+1}|^2
=1.
\]
Similarly,
\begin{equation}\label{eq:cmv-factor-even}
R_{2k+1}(z)S_{2k+2}(z)
=
W_{2k+2}(z)R_{2k+2}(z),
\end{equation}
where
\[
W_{2k+2}(z)
=
\begin{pmatrix}
0&\rho_{2k+1}^2z^2\\
1&
\overline{\alpha_{2k+2}}
+
\overline{\alpha_{2k+1}}z
\end{pmatrix}
\]
and
\[
R_{2k+2}(z)
=
\begin{pmatrix}
0&\rho_{2k+2}^2\\
z&\alpha_{2k+2}
\end{pmatrix}.
\]
Thus \(R_{2k+2}\) has exactly the same form as \(R_{2k}\), and the factorization repeats with period two.

Consequently,
\[
H(z)S_0(z)S_1(z)\cdots S_N(z)
=
W_0(z)W_1(z)\cdots W_N(z)R_N(z).
\]
Hence, apart from the terminal remainder, the Schur product has been converted into a product of continued-fraction matrices. This yields the following formal continued-fraction representation of the Schur algorithm.

\begin{theorem}\label{thm:cmv-cf}
Let \(f\) be a Schur function with Schur parameters $\{\alpha_k\}_{k=0}^{\infty}$ and as before let
\[
\rho_n^2=1-|\alpha_n|^2.
\]
Then the Schur algorithm generates the formal continued fraction
\[
\frac{f(z)}{1-zf(z)}
=
\cfrac{1}{
\overline{\alpha_0}-z+
\cfrac{\rho_0^2}{
\alpha_0+\alpha_1z+
\cfrac{\rho_1^2z^2}{
\overline{\alpha_2}+\overline{\alpha_1}z+
\cfrac{\rho_2^2}{
\alpha_2+\alpha_3z+
\cfrac{\rho_3^2z^2}{
\overline{\alpha_4}+\overline{\alpha_3}z+\ddots
}}}}}.
\]
Here the equality is understood formally: for every finite number of Schur steps, the corresponding finite matrix product yields the associated finite continued fraction, with the remaining Schur iterate appearing as the terminal parameter.
\end{theorem}

As usual, the continued fraction generates a three-term recurrence. Let \(Q_n(z)\) denote its denominator continuants, normalized by
\[
Q_{-1}(z)=0,\qquad Q_0(z)=1,
\qquad
Q_1(z)=\overline{\alpha_0}-z.
\]
Then
\[
Q_{2k+2}(z)
=
\bigl(\alpha_{2k}+\alpha_{2k+1}z\bigr)Q_{2k+1}(z)
+
\rho_{2k}^2Q_{2k}(z),
\qquad k\geq0,
\]
while
\[
Q_{2k+3}(z)
=
\bigl(\overline{\alpha_{2k+2}}
+\overline{\alpha_{2k+1}}z\bigr)Q_{2k+2}(z)
+
\rho_{2k+1}^2z^2Q_{2k+1}(z),
\qquad k\geq0.
\]
Thus the recurrence alternates between constant and quadratic coefficients,
\[
\rho_{2k}^2,
\qquad
\rho_{2k+1}^2z^2,
\]
reflecting the interpolation pattern
\[
\infty,\infty,0,0,\infty,\infty,0,0,\ldots .
\]
Equivalently, the recurrence is encoded by the tridiagonal matrix
\[
T(z)=
\begin{pmatrix}
\overline{\alpha_0}-z
&\rho_0&&&&\\
-\rho_0
&\alpha_0+\alpha_1z
&z\rho_1&&&\\
&-z\rho_1
&\overline{\alpha_2}+\overline{\alpha_1}z
&\rho_2&&\\
&&-\rho_2
&\alpha_2+\alpha_3z
&z\rho_3&\\
&&&-z\rho_3
&\overline{\alpha_4}+\overline{\alpha_3}z
&\ddots\\
&&&&\ddots&\ddots
\end{pmatrix}.
\]
Indeed, the leading principal minors of \(T(z)\) satisfy precisely the two recurrence relations above.

The CMV structure becomes visible after a simple normalization. Let
\[
D=\operatorname{diag}(1,-1,1,-1,\ldots)
\]
and introduce the unitary \(2\times2\) matrices
\[
\Theta(\alpha_n)
=
\begin{pmatrix}
\overline{\alpha_n}&\rho_n\\
\rho_n&-\alpha_n
\end{pmatrix},
\qquad
\rho_n=\sqrt{1-|\alpha_n|^2}.
\]
Define
\[
\mathcal L
=
\operatorname{diag}(\Theta(\alpha_0),
\Theta(\alpha_2), \Theta(\alpha_4)\cdots)
\]
and
\[
\mathcal M
=
\operatorname{diag}(1,
\Theta(\alpha_1),
\Theta(\alpha_3),\cdots ).
\]
A direct comparison of the entries gives
\begin{equation}\label{eq:cmv-pencil-factor}
DT(z)=\mathcal L-z\mathcal M^*.
\end{equation}
Thus the three-term recurrence generated by the continued fraction is, up to the harmless diagonal factor \(D\), the generalized eigenvalue problem
\begin{equation}\label{eq:cmv-gevp}
\mathcal L u=z\mathcal M^*u.
\end{equation}
Pencils of this form appeared earlier in the numerical linear algebra literature under the name \emph{Schur parameter pencils}. They were introduced by Bunse-Gerstner and Elsner in connection with the unitary eigenvalue problem \cite{BGE91}. The matrices \(\mathcal L\) and \(\mathcal M\) are precisely the two block-diagonal factors occurring in the standard CMV factorization
\[
\mathcal C=\mathcal L\mathcal M.
\]
For the CMV representation and its relation to Schur parameters and orthogonal polynomials on the unit circle, see \cite{SimonOPUC}.

To complete the derivation, note that, since \(\mathcal M\) is unitary, multiplying \eqref{eq:cmv-gevp} by \(\mathcal M\) gives
\[
\mathcal M\mathcal L\,u=zu.
\]
The matrix
\[
\widetilde{\mathcal C}=\mathcal M\mathcal L
\]
is the alternate CMV matrix. Thus the continued fraction constructed above leads directly, through its denominator recurrence, to a CMV-type generalized eigenvalue problem and hence to a CMV matrix. With the opposite Laurent ordering one obtains instead the standard CMV matrix
\[
\mathcal C=\mathcal L\mathcal M.
\]

\section{Nevanlinna--Pick algorithm and \(R_{II}\)-fractions}

We are now ready to bring the Thiele construction and the Schur algorithm
together. In the form discussed above, Wall's construction already does
this when the Schur transformation is repeatedly applied at a fixed
interpolation point. To accommodate a sequence of distinct interpolation
nodes, we turn to the Nevanlinna--Pick interpolation problem.

In \cite{DZ09}, the modified Nevanlinna--Pick transformation was derived
directly from the Schwarz lemma, without invoking Wall's construction.
This is essentially the mechanism underlying the classical Schur
algorithm, but the interpolation point may now vary from step to step.
Iterating the resulting transformations leads to an $R_{II}$-fraction.

Recall that a Nevanlinna function is analytic in $\mathbb C_+$ and maps
$\mathbb C_+$ into itself. We extend it to $\mathbb C_-$ by symmetry. Also, in what follows, we assume that a Nevanlinna function is not rational.
Let $\varphi$ be such a function and apply the classical Nevanlinna--Pick
algorithm to $f_0=-1/\varphi$ at the interpolation points
\[
z_j=\xi_j+iy_j\in\mathbb C_+.
\]
Write $f_j(z_j)=u_j+iv_j$, and assume that the steps under consideration
are nondegenerate, so that $v_j>0$. Put
\[
\ell_j(\lambda)=\frac{\lambda-\xi_j}{y_j},\qquad
\alpha_j(\lambda)=u_j+v_j\ell_j(\lambda),
\]
and
\[
\beta_j=\frac{\sqrt{v_j}}{y_j}>0,\qquad
q_j(\lambda)=(\lambda-z_j)(\lambda-\overline{z_j}).
\]
The classical Nevanlinna-Pick transformation \cite{Ach61} then reads
\begin{equation}\label{eq:NP-classical-step}
f_j(\lambda)=u_j+
v_j\frac{\ell_j(\lambda)f_{j+1}(\lambda)-1}
{f_{j+1}(\lambda)+\ell_j(\lambda)}.
\end{equation}
The reduced function has nonnegative imaginary part in $\mathbb C_+$,
as follows from the Schwarz lemma. Real constant and infinite terminal
parameters are interpreted in the usual way.

For $M=\begin{pmatrix}a&b\\c&d\end{pmatrix}$, write
$\mathcal T_M[w]=(aw+b)/(cw+d)$. Passing to $h_j=-1/f_j$, with
$h_0=\varphi$, gives
\[
h_j=\mathcal T_{M_j}[h_{j+1}],\qquad
M_j(\lambda)=
\begin{pmatrix}
\ell_j(\lambda)&-1\\
v_j-u_j\ell_j(\lambda)&\alpha_j(\lambda)
\end{pmatrix}.
\]

\begin{proposition}\label{prop:NP-factorization}
For each $j$, we have
\[
M_j(\lambda)=C_j(\lambda)L_j(\lambda),
\]
where
\[
C_j(\lambda)=
\begin{pmatrix}
0&-1\\
\beta_j^2q_j(\lambda)&\alpha_j(\lambda)
\end{pmatrix},
\qquad
L_j(\lambda)=
\begin{pmatrix}
1&0\\
-\ell_j(\lambda)&1
\end{pmatrix}.
\]
Consequently,
\[
L_j(\lambda)C_{j+1}(\lambda)=
\begin{pmatrix}
0&-1\\
\beta_{j+1}^2q_{j+1}(\lambda)&
\alpha_{j+1}(\lambda)+\ell_j(\lambda)
\end{pmatrix}.
\]
In particular, $L_jC_{j+1}$ is again a continued-fraction matrix,
with an affine partial denominator and a quadratic partial numerator.
\end{proposition}

\begin{proof}
Since $\beta_j^2q_j=v_j(1+\ell_j^2)$, both identities follow by direct
multiplication.
\end{proof}

The proposition gives
\[
M_0\cdots M_{N-1}
=
C_0(L_0C_1)\cdots(L_{N-2}C_{N-1})L_{N-1}.
\]
All factors except the last are continued-fraction matrices of
$R_{II}$ type. The last factor replaces $h_N$ by the terminal parameter
\[
\tau_N=\mathcal T_{L_{N-1}}[h_N]
=-\frac1{f_N+\ell_{N-1}},
\]
so the finite representation is exact with this terminal tail.
Here $\tau_N$ is a Nevanlinna function when $f_N$ is finite, and
$\tau_N=0$ when $f_N\equiv\infty$.

The preceding calculation uses one normalization of the classical
parameters. For the matrix realization, we now use the
modified transformation of \cite{DZ09,D10}. Suppose that $\varphi_0=\varphi$
is a Cauchy transform of a probability measure on $\mathbb R$. At each
step, put
\begin{equation}\label{h1forS}
\Phi_j(\lambda)=-\frac1{\varphi_j(\lambda)},\qquad
w_j=\Phi_j(z_j),
\end{equation}
and define
\begin{equation}\label{h2forS}
a_j^{(2)}=\frac{\Im w_j}{\Im z_j},\qquad
a_j^{(1)}=a_j^{(2)}z_j-w_j,\qquad
b_j^2=a_j^{(2)}-1.
\end{equation}
At a nonterminating step, $a_j^{(2)}>1$ and we choose $b_j>0$.
With this choice, the reduced function in
\begin{equation}\label{eq:NP-normalized-step}
\varphi_j(\lambda)=
-\frac1{a_j^{(2)}\lambda-a_j^{(1)}
+b_j^2q_j(\lambda)\varphi_{j+1}(\lambda)}
\end{equation}
is again a Cauchy transform of a probability measure and so it is a Nevanlinna function. We henceforth use $a_j^{(1)},a_j^{(2)},b_j$ for
these normalized coefficients, distinguishing them from the classical
parameters above.

Writing
\[
d_j(\lambda)=a_j^{(2)}\lambda-a_j^{(1)},\qquad
W_j(\lambda)=
\begin{pmatrix}
0&-1\\
b_j^2q_j(\lambda)&d_j(\lambda)
\end{pmatrix},
\]
we obtain
\begin{equation}\label{eq:NP-product}
\varphi=\mathcal T_{W_0\cdots W_{N-1}}[\varphi_N].
\end{equation}
Thus the modified Nevanlinna--Pick transformation has the same
continued-fraction matrix structure as the preceding refactorization.
As in Wall's normalization, the partial numerators contain conjugate
pairs of interpolation points; here the fixed pair $i,-i$ is replaced
by $z_j,\overline{z_j}$. The resulting formal continued fraction is
\begin{equation}\label{eq:NP-RII}
-\cfrac{1}{
a_0^{(2)}\lambda-a_0^{(1)}
-\cfrac{b_0^2(\lambda-z_0)(\lambda-\overline z_0)}
{a_1^{(2)}\lambda-a_1^{(1)}
-\cfrac{b_1^2(\lambda-z_1)(\lambda-\overline z_1)}
{a_2^{(2)}\lambda-a_2^{(1)}-\ddots}}}.
\end{equation}
At each finite stage, the identity is \eqref{eq:NP-product}, with the
actual terminal tail $\varphi_N$.

As before, the continued fraction \eqref{eq:NP-RII} generates a
three-term recurrence relation. Let
\[
P_0(\lambda)=1,\qquad
P_1(\lambda)=a_0^{(2)}\lambda-a_0^{(1)}.
\]
Then the denominator polynomials of its convergents satisfy
\begin{equation}\label{eq:NP-RII-rec}
P_{j+1}(\lambda)=
\bigl(a_j^{(2)}\lambda-a_j^{(1)}\bigr)P_j(\lambda)
-b_{j-1}^2(\lambda-z_{j-1})(\lambda-\overline z_{j-1})
P_{j-1}(\lambda),
\quad j\geq1.
\end{equation}
For nonterminating steps, introduce
\[
\widehat P_j(\lambda)=
\frac{P_j(\lambda)}
{b_0\cdots b_{j-1}\prod_{k=0}^{j-1}(z_k-\lambda)},
\qquad j\geq1,\qquad \widehat P_0(\lambda)=1.
\]
Then \eqref{eq:NP-RII-rec} becomes
\begin{equation}\label{eq:NP-symmetric-rec}
b_j(\lambda-z_j)\widehat P_{j+1}(\lambda)
+\bigl(a_j^{(2)}\lambda-a_j^{(1)}\bigr)\widehat P_j(\lambda)
+b_{j-1}(\lambda-\overline z_{j-1})\widehat P_{j-1}(\lambda)=0,
\quad j\geq1.
\end{equation}
For $j=0$, the relation is
\[
b_0(\lambda-z_0)\widehat P_1(\lambda)
+\bigl(a_0^{(2)}\lambda-a_0^{(1)}\bigr)\widehat P_0(\lambda)=0.
\]

Consequently, with
$\widehat P=(\widehat P_0,\widehat P_1,\ldots)^\top$, the recurrence
can be written as the formal generalized eigenvalue equation
\[
J^{(1)}\widehat P=\lambda J^{(2)}\widehat P,
\]
where
\[
J^{(1)}=
\begin{pmatrix}
a_0^{(1)}&b_0z_0&&\\
b_0\overline z_0&a_1^{(1)}&b_1z_1&\\
&b_1\overline z_1&a_2^{(1)}&\ddots\\
&&\ddots&\ddots
\end{pmatrix}.
\]
Since $a_j^{(2)}=1+b_j^2$, the second matrix has the factorization
\begin{equation}\label{eq:NP-J2-Cholesky}
J^{(2)}=
\begin{pmatrix}
a_0^{(2)}&b_0&&\\
b_0&a_1^{(2)}&b_1&\\
&b_1&a_2^{(2)}&\ddots\\
&&\ddots&\ddots
\end{pmatrix}
=
\underbrace{\begin{pmatrix}
1&b_0&&\\
&1&b_1&\\
&&1&\ddots\\
&&&\ddots
\end{pmatrix}}_{\mathcal L^*}
\underbrace{\begin{pmatrix}
1&&&\\
b_0&1&&\\
&b_1&1&\\
&&\ddots&\ddots
\end{pmatrix}}_{\mathcal L}.
\end{equation}
Both matrices are Hermitian. The factorization also shows that every
finite leading principal submatrix of $J^{(2)}$ is positive definite.
Thus the Schur transformation leads, through the same continued-fraction
mechanism as before, to a pencil of tridiagonal Hermitian matrices.

We now assume that the representing probability measure has infinite
support contained in a compact interval $[\alpha,\beta]$, so that
\[
\varphi(\lambda)=\int_\alpha^\beta\frac{d\sigma(t)}{t-\lambda}.
\]
Let $J^{(1)},J^{(2)}$ be the matrices obtained from the
modified Nevanlinna--Pick transformation, and let
$e_0=(1,0,\ldots)^\top\in\ell^2(\mathbb N_0)$.

\begin{proposition}\label{prop:NP-mfunction}
Assume that $z_j\to\infty$. Then $J^{(1)}$ and $J^{(2)}$ define
bounded self-adjoint operators in $\ell^2(\mathbb N_0)$, $J^{(2)}\geq cI$ for some $c>0$, and
\begin{equation}\label{eq:NP-mfunction}
\varphi(\lambda)=
\left\langle
\bigl(J^{(1)}-\lambda J^{(2)}\bigr)^{-1}e_0,e_0
\right\rangle_{\ell^2(\mathbb N_0)}=\int_\alpha^\beta\frac{d\sigma(t)}{t-\lambda},
\qquad \lambda\in\mathbb C\setminus[\alpha,\beta].
\end{equation}
In particular, $\varphi$ is the $m$-function of the pencil
$J^{(1)}-\lambda J^{(2)}$.
\end{proposition}

\begin{proof}
By \cite[Proposition~5.2 and Theorem~5.4]{BDZh}, the matrices define
bounded self-adjoint operators and $0\notin\sigma(J^{(2)})$.
Positivity of $J^{(2)}$ therefore gives $J^{(2)}\geq cI$ for some $c>0$.
Consequently,
\[
A=(J^{(2)})^{-1/2}J^{(1)}(J^{(2)})^{-1/2}
\]
is bounded and self-adjoint, and
\[
\bigl(J^{(1)}-\lambda J^{(2)}\bigr)^{-1}
=(J^{(2)})^{-1/2}(A-\lambda I)^{-1}(J^{(2)})^{-1/2}.
\]
Denote the right-hand side of \eqref{eq:NP-mfunction} by $m(\lambda)$.
Schur complements reproduce \eqref{eq:NP-normalized-step}, so
$m(z_j)=\varphi(z_j)$. Boundedness of $J^{(1)}$ and invertibility of
$J^{(2)}$ imply that $m$ is analytic at infinity. The same holds for
$\varphi$ because $\sigma$ has compact support. Since the interpolation
points tend to infinity, the identity theorem gives $m=\varphi$ near
infinity and then throughout $\mathbb C\setminus[\alpha,\beta]$.
\end{proof}

The preceding construction also gives a Stieltjes-type continued
fraction. In this specialization, we write $a_j=a_j^{(2)}$.

\begin{theorem}\label{thm:Stieltjes-CF}
Let
\[
s(\lambda)=\int_{[0,\infty)}\frac{d\sigma(t)}{t-\lambda},
\]
where $\sigma$ is a probability measure, and let $\{x_k\}_{0}^{\infty}$
be a sequence of distinct positive numbers. Then the normalized reduction generates
the formal continued fraction
\begin{equation}\label{eq:Stieltjes-CF}
s(\lambda)=
-\cfrac{1}{
a_0\lambda
-\cfrac{b_0^2(\lambda+x_0)}
{a_1
-\cfrac{b_1^2(\lambda+x_1)}
{a_2\lambda
-\cfrac{b_2^2(\lambda+x_2)}
{a_3-\ddots}}}}.
\end{equation}
At every nonterminating step,
\[
a_j>1,\qquad b_j>0,\qquad b_j^2=a_j-1.
\]
If $a_j=1$, then $b_j=0$ and the fraction terminates.
The formal equality means that every finite sequence of reductions
gives an exact finite continued fraction with the remaining transformed
function as its terminal parameter.
\end{theorem}

Thus symmetry reduces the $R_{II}$-fraction to an $R_I$-type fraction
with alternating linear and constant partial denominators.

\begin{proof}
The function
\[
\varphi(\zeta)=\zeta s(\zeta^2)
=\frac12\int_{[0,\infty)}
\left(\frac1{\sqrt t-\zeta}+\frac1{-\sqrt t-\zeta}\right)d\sigma(t)
\]
is an odd Nevanlinna function and is the Cauchy transform of a symmetric
probability measure on $\mathbb R$.
Apply the normalized reduction at $z_j=i\sqrt{x_j}$.
If $\varphi_j$ is odd, then $\Phi_j=-1/\varphi_j$ is also odd.
Together with symmetry, this implies that
$\Phi_j(i\sqrt{x_j})$ is purely imaginary. Hence \eqref{h2forS} gives
$a_j^{(1)}=0$. At a nonterminating step,
\eqref{eq:NP-normalized-step} becomes
\[
\varphi_j(\zeta)=
-\frac1{a_j\zeta+b_j^2(\zeta^2+x_j)\varphi_{j+1}(\zeta)},
\qquad b_j^2=a_j-1.
\]
Equivalently,
\[
\varphi_{j+1}(\zeta)=
\frac{-1/\varphi_j(\zeta)-a_j\zeta}
{b_j^2(\zeta^2+x_j)}.
\]
This identity preserves oddness, so every nonterminal iterate is odd.
Write
\[
\varphi_j(\zeta)=\zeta s_j(\zeta^2),\qquad s_0=s,
\]
where each $s_j$ is a normalized Stieltjes function, and define
\[
r_j(\lambda)=
\begin{cases}
s_j(\lambda),&j\ \text{even},\\
\lambda s_j(\lambda),&j\ \text{odd},
\end{cases}
\qquad
d_j(\lambda)=
\begin{cases}
a_j\lambda,&j\ \text{even},\\
a_j,&j\ \text{odd}.
\end{cases}
\]
Substituting $\lambda=\zeta^2$ in the reduction gives
\[
r_j(\lambda)=
-\frac1{d_j(\lambda)+b_j^2(\lambda+x_j)r_{j+1}(\lambda)}.
\]
Iterating this identity yields \eqref{eq:Stieltjes-CF}, with $r_N$
as the terminal parameter after $N$ steps.
\end{proof}

\begin{remark}
If all moments of $\sigma$ are finite, the corresponding reductions
at infinity produce the classical Stieltjes continued fraction.
Here interpolation at infinity is understood through the appropriately
rescaled limiting reduction, rather than by direct substitution of
$x_j=\infty$ in \eqref{eq:Stieltjes-CF}.
\end{remark}

The continued fraction \eqref{eq:Stieltjes-CF} also gives rise to a
tridiagonal linear pencil. Let $P_j$ now denote its denominator
polynomials, initialized by
\[
P_0(\lambda)=1,\qquad P_1(\lambda)=d_0(\lambda)=a_0\lambda.
\]
With $d_j$ as in the proof, they satisfy
\[
P_{j+1}(\lambda)=d_j(\lambda)P_j(\lambda)
-b_{j-1}^2(\lambda+x_{j-1})P_{j-1}(\lambda),
\qquad j\geq1.
\]
On the nondegenerate part of this recurrence, set
\[
\widehat P_j(\lambda)=
\frac{P_j(\lambda)}{b_0\cdots b_{j-1}},
\qquad \widehat P_0(\lambda)=1.
\]
Then
\[
b_j\widehat P_{j+1}(\lambda)
-d_j(\lambda)\widehat P_j(\lambda)
+b_{j-1}(\lambda+x_{j-1})\widehat P_{j-1}(\lambda)=0,
\qquad j\geq1,
\]
with boundary relation
$b_0\widehat P_1-d_0\widehat P_0=0$.
Consequently, the formal equation
$\mathcal T(\lambda)\widehat P=0$ uses the pencil
\[
\mathcal T(\lambda)=
\begin{pmatrix}
a_0\lambda&b_0&&\\
b_0(\lambda+x_0)&a_1&b_1&\\
&b_1(\lambda+x_1)&a_2\lambda&\ddots\\
&&\ddots&\ddots
\end{pmatrix}.
\]
The alternating dependence on $\lambda$ along the diagonal comes from
the parity reduction $\varphi(\zeta)=\zeta s(\zeta^2)$.
In this normalization, the coefficient matrices of $\mathcal T(\lambda)$
are generally not Hermitian.

\section{A final detour: periodic nodes and banded matrices}

Let us return once more to the pencil
\[
J^{(1)}-\lambda J^{(2)}
\]
associated with the multipoint Nevanlinna--Pick algorithm, now allowing
interpolation nodes at infinity. We assume that the reductions under
consideration are nondegenerate and that the moments required for the
steps at infinity are finite. For a finite interpolation point $z_j$,
one step of the transformation in question has the form
\[
\varphi_j(\lambda)=
-\frac{1}{
a_j^{(2)}\lambda-a_j^{(1)}
+b_j^2(\lambda-z_j)(\lambda-\overline{z_j})
\varphi_{j+1}(\lambda)}.
\]
At the node $z_j=\infty$, the transformation is known and becomes the classical Schur transformation in the Nevanlinna class 
\[
\varphi_j(\lambda)=
-\frac{1}{
\lambda-a_j+\widehat b_j^{\,2}\varphi_{j+1}(\lambda)},
\qquad a_j\in\mathbb R,\quad \widehat b_j>0
\]
which generates a Jacobi fraction provided that the moments of the underlying measure are finite.

This observation is particularly transparent at the level of the
associated recurrence. For a finite node, the off-diagonal terms depend
on $\lambda$ through
\[
b_j(\lambda-z_j),
\qquad
b_j(\lambda-\overline{z_j}),
\]
whereas at an infinite node they are constant. Consequently, when
$z_j=\infty$, the corresponding entries of the pencil satisfy
\[
(J^{(2)})_{j,j+1}=(J^{(2)})_{j+1,j}=0,
\]
while the corresponding off-diagonal entries of $J^{(1)}$ remain
nonzero.

Suppose now that the node at infinity occurs at regularly spaced
indices:
\[
z_{p-1}=z_{2p-1}=z_{3p-1}=\cdots=\infty.
\]
Then $J^{(2)}$ decomposes into a direct sum of $p\times p$ tridiagonal
blocks,
\[
J^{(2)}=\operatorname{diag}(B_0,B_1,B_2,\ldots).
\]
Under the normalization adopted above, each block $B_k$ is positive
definite. Hence its blockwise inverse is
\[
\bigl(J^{(2)}\bigr)^{-1}
=\operatorname{diag}(B_0^{-1},B_1^{-1},B_2^{-1},\ldots).
\]
This inverse defines a bounded operator on $\ell^2(\mathbb N_0)$
precisely when
\[
\inf_{k\geq0}\lambda_{\min}(B_k)>0.
\]
The matrix identities below are understood blockwise and do not
require this uniform bound.

Take the blockwise Cholesky factorizations
\[
J^{(2)}=LL^*=UU^*,
\]
where
\[
L=\operatorname{diag}(L_0,L_1,\ldots),
\qquad
U=\operatorname{diag}(U_0,U_1,\ldots),
\]
each $L_k$ is lower triangular, and each $U_k$ is upper triangular.
The second factorization is the reverse Cholesky factorization.
We write $L^{-*}=(L^*)^{-1}$, and similarly for $U$.

\begin{theorem}\label{thm:RII-banded}
Assume that $z_{kp+p-1}=\infty$ for every $k\geq0$. Then
\[
A_L=L^{-1}J^{(1)}L^{-*},
\qquad
A_U=U^{-1}J^{(1)}U^{-*}
\]
are Hermitian matrices of bandwidth at most $p$; that is, they are
$(2p+1)$-diagonal. Moreover, the formal generalized eigenvalue problem
\[
J^{(1)}u=\lambda J^{(2)}u
\]
is equivalent to
\[
A_Lv=\lambda v,\qquad v=L^*u,
\]
and also to
\[
A_Uw=\lambda w,\qquad w=U^*u.
\]
\end{theorem}

\begin{proof}
Both transformed matrices are Hermitian, and the stated changes of
variables give the equivalence of the eigenvalue equations.

Partition $J^{(1)}$ into $p\times p$ blocks. Its only nonzero blocks
outside the block diagonal join consecutive blocks and have the form
\[
c_k e_{p-1}e_0^*,
\]
where $e_0,\ldots,e_{p-1}$ are the standard coordinate vectors in
$\mathbb C^p$. Thus the block of $A_L$ joining blocks $k$ and $k+1$ is
\[
c_k(L_k^{-1}e_{p-1})(L_{k+1}^{-1}e_0)^*.
\]
Since $L_k^{-1}$ is lower triangular, $L_k^{-1}e_{p-1}$ is a scalar
multiple of $e_{p-1}$. This block is therefore supported in its last
row, so its entries lie within the first $p$ upper diagonals.
The diagonal blocks lie within bandwidth $p-1$, and Hermitian symmetry
gives the same bound below the diagonal. Hence
\[
(A_L)_{jk}=0\qquad\text{whenever}\qquad |j-k|>p.
\]

For the reverse Cholesky factorization, $U_{k+1}^{-1}e_0$ is a scalar
multiple of $e_0$. The corresponding off-diagonal block of $A_U$ is
therefore supported in its first column, which gives the same
bandwidth bound.
\end{proof}

\begin{remark}
One may also work directly with
\[
A=\bigl(J^{(2)}\bigr)^{-1}J^{(1)}.
\]
Then the original generalized eigenvalue equation is equivalent to
$Au=\lambda u$. This matrix also satisfies
\[
A_{jk}=0\qquad\text{whenever}\qquad |j-k|>p,
\]
because left multiplication by the block-diagonal inverse can spread
each interblock coupling over only one block. Although $A$ is
generally not Hermitian, it is quasiseparable, for example, see \cite{Vandebril08}.
\end{remark}

Thus the regularly spaced occurrence of the interpolation node at
infinity converts an $R_{II}$-pencil into an ordinary matrix with
finite bandwidth. The reduction is not canonical. One may use
$\bigl(J^{(2)}\bigr)^{-1}J^{(1)}$ directly or preserve the Hermitian
structure through either Cholesky reduction. Other blockwise
factorizations give further banded representatives of the same
formal generalized eigenvalue problem. This freedom is one reason
why we do not attempt an entry-by-entry identification with a
particular banded-matrix model.

For $p=1$, every interpolation node is at infinity, $J^{(2)}=I$,
and one recovers the usual Jacobi matrix. More generally, if
\[
p=g+1,
\]
then the resulting matrices are $(2g+3)$-diagonal. This is reminiscent
of banded-matrix constructions arising in the theory of orthogonal
rational functions and Laurent polynomials; see, in particular,
the work of Simonov \cite{Simonov06,Simonov06Moment} and the GMP
matrices introduced by Yuditskii \cite{Yuditskii18}. GMP matrices
associated with $g$ prescribed finite real poles and the point at
infinity have the same bandwidth. An identification with a
particular model requires matching its pole data, normalization,
basis, and defining structural conditions. The observation above
exhibits a mechanism behind the finite-band structure: periodic
returns to the node at infinity split $J^{(2)}$ into finite blocks
and permit a banded realization of the pencil.

\section*{Epilogue}

A more systematic treatment of the spectral theory may follow, if anyone cares enough to ask.

\section*{Tool and computational resource disclosure}
During the preparation of this manuscript, the author used ChatGPT (OpenAI) to explore alternative derivations, check algebraic manipulations, identify potentially relevant literature, and improve the exposition. The author independently verified all AI-assisted mathematical arguments and bibliographic information included in the final manuscript and takes full responsibility for the paper’s content and correctness.

\end{document}